\documentclass[12pt]{amsart}
\usepackage{amsthm,amsfonts, amssymb, amscd}

\newtheorem{theorem}{\bf Theorem}
\newtheorem{lemma}[theorem]{\bf Lemma}
\newtheorem{propos}[theorem]{\bf Proposition}

\theoremstyle{definition}

\newtheorem{rmk}[theorem]{\bf Remark}

\numberwithin{subsection}{section}

\newcommand{\C}{\mathbb C}

\newcommand{\Pro}{\mathbb P}
\newcommand{\Q}{\mathbb Q}
\newcommand{\Z}{\mathbb Z}

\newcommand{\Oc}{\mathcal O}
\newcommand{\Qc}{\mathcal Q}
\newcommand{\Cf}{\mathfrak C}
\newcommand{\Mf}{\mathfrak M}
\newcommand{\X}{\mathfrak X}

\newcommand{\Gm}{\mathbb G_{\textbf{m}}}

\newcommand{\sym}{\text{Sym}}

\DeclareMathOperator{\GL}{GL}
\DeclareMathOperator{\PGL}{PGL}
\DeclareMathOperator{\Spec}{Spec}
\DeclareMathOperator{\Proj}{Proj}
\DeclareMathOperator{\Sym}{Sym}
\DeclareMathOperator{\Image}{Im}
\DeclareMathOperator{\Hom}{Hom}

\begin{document}

\title[The integral Chow ring of $\Mf^{\leq 1}_0$.]{The integral Chow
ring of the stack of at most 1-nodal rational curves}

\author{Dan Edidin and Damiano Fulghesu}
\address{Department of Mathematics, University of Missouri, Columbia, MO 65211}
\email{edidin@math.missouri.edu, damiano@math.missouri.edu}

\begin{abstract} We give a presentation for the stack $\Mf^{\leq 1}_0$
of rational curves with at most 1 node as the quotient by $\GL_3$ of
an open set in a 6 dimensional irreducible representation.  We then
use equivariant intersection theory to calculate the integral Chow
ring of this stack.
\end{abstract}
\maketitle

\section{Introduction}
The integral Chow ring of a smooth quotient stack was defined in the
paper \cite{eit}. Subsequently, Kresch \cite{kr} extended the
definition to smooth stacks which admit stratifications by quotient
stacks.  In \cite{eit} the integral Chow rings of the stacks of
elliptic curves $\Mf_{1,1}$ and $\overline{\Mf}_{1,1}$ were
computed. In \cite{vi-ap} Vistoli calculated the integral Chow ring of
$\Mf_2$. All of these calculations use presentations for the stacks as
quotient stacks.

In this paper we turn our attention to stacks of unpointed rational
curves. Unlike stacks of curves of positive genus, they are not
Deligne-Mumford, since $\Pro^1$ has an infinite automorphism group. In
\cite{fulg} the second author introduced the stratification of the
stack of rational\footnote{that is to say of arithmetic genus 0}
curves $\Mf_0$ by the number of nodes.  The stack $\Mf_0^0$ of smooth
rational curves is the classifying stack $B\PGL_2$ whose Chow ring was
computed by Pandharipande \cite{Pan}. The focus of this paper is
$\Mf_0^{\leq 1}$, the stack of rational curves with at most one node.

Although it is not a Deligne-Mumford stack, $\Mf_0^{\leq 1}$ 
admits a relatively
easy description as a quotient. On a 
rational curve $C$ with at most one node the dual of the dualizing
sheaf embeds the curve as a plane conic.
This observation
allows us to prove that 
(Proposition \ref{prop.presentation}) that $\Mf_0^{\leq 1}$
is a quotient by $\GL_3$ of an open set in the representation $
H^0(\omega^{\vee}_{\Pro^2}(-1))$. We can then calculate the Chow ring
of the stack as the $\GL_3$-equivariant Chow ring of this open
set and obtain the following result.

Let $\Cf \stackrel{\pi} \to \Mf^{\leq 1}_0$ be the universal
curve and let $E = \pi_* \omega_\pi^\vee$.
The vector bundle $E$ is locally free of rank $3$,
so we denote the Chern classes of $E$ by $c_1, c_2 ,c_3$ respectively.
\begin{theorem}
$$A^* \Mf^{\leq 1}_0 = \Z[c_1, c_2, c_3]/(4c_3, 2 c_1 c_3, c_1^2 c_3)$$
\end{theorem}

For any algebraic space 
the $\GL_3$-equivariant Chow ring injects into the $T$-equivariant Chow ring,
where $T$ is a maximal torus. Thus we may check relations in $T$-equivariant
Chow ring. This allows us to use the localization theorem for torus
actions and reduce
the problem to one of polynomial interpolation.
In the last section of the paper we use the techniques developed
for the calculation of $A^*(\Mf^{\leq 1}_0)$ to compute Chow rings
of stacks of quadrics in $\Pro^{n-1}$ and classifying spaces
for certain extensions of $SO(n)$ by a finite group.

\begin{rmk}
If a rational curve $C$ has two or more 
nodes then $\omega_C^{\vee}$ is no longer ample and there
is no obvious presentation for the stack $\Mf_0^{\leq n}$
of rational curves with at most $n$ nodes as a quotient stack.
In fact there is some evidence that
$\Mf_0^{\leq n}$ {\em is not} a quotient stack since the second author proved
\cite{fulg} that if
$n \geq 2$ the projection from the universal curve $\Cf_0^{\leq n} \to
\Mf_0^{\leq n}$ is not representable in the category of schemes.
However, using the intersection theory developed by Kresch for stacks
stratified by quotient stacks, the second author \cite{fulg} was able to compute
the {\em rational} Chow rings of the stacks
$\Mf_0^{\leq 2}$ and $\Mf_0^{\leq 3}$, but it is unclear how to 
determine their integral Chow rings.
\end{rmk}

\section{Some facts about equivariant Chow groups}
Equivariant Chow groups were defined in the paper \cite{eit}.
We briefly 
recall some basic facts and notation that we will use in our computation.

We work over an arbitrary field $K$.
Let $G$ be a linear algebraic group. 
For any algebraic space $X$ we denote
the direct sum of the equivariant Chow groups by $A_*^G(X)$.
If $X$ is smooth then there is a product structure on
equivariant Chow groups and we denote the equivariant Chow
ring by the notation $A^*_G(X)$. 
Following standard notation we denote the equivariant Chow ring
of a point by $A^*_G$.
Flat pullback 
via the morphism $X \to \Spec K$ 
makes the equivariant Chow groups $A_*^G(X)$ into
an $A^*_G$-module. When $X$ is smooth the equivariant Chow ring 
$A^*_G(X)$ becomes an $A^*_G$-algebra.

The relation between equivariant Chow rings and Chow rings
of quotient stacks is given by the following result.
\begin{propos}\cite[Propositions 17, 19]{eit}
If ${\mathcal F} = [X/G]$ then the equivariant Chow ring $A^*_G(X)$
is independent of the presentation for ${\mathcal F}$ and may
be identified with the operational Chow ring of ${\mathcal F}$.
\end{propos}
If $V$ is a $G$-module, then $V$ defines a $G$-equivariant vector bundle
over $\Spec K$. Consequently a representation $V$ of rank
$r$ has {\it Chern classes} $c_1, \ldots , c_r \in A^*_G$. If
$X$ is a smooth algebraic space then we will view the Chern
classes as elements of the equivariant Chow ring $A^*_G(X)$ via 
the pullback $A^*_G \to A^*_G(X)$.

\begin{rmk} For the presentation 
of the stack $\Mf_0^{\leq 1}$ we work over $\Spec \Z$. For
the calculation of Chow rings we work over an arbitrary
field.
\end{rmk}

\subsection{Equivariant Chow rings for $\GL_n$ actions}
Let $T = \Gm^n$ be a maximal torus.  Because $\GL_n$ is {\em
special}\footnote{This means that every $\GL_n$-torsor is locally
trivial in the Zariski topology.} the restriction homomorphism
$A^*_{\GL_n} \to A^*_T$ is injective and the image consists of the
classes invariant under the action of the Weyl group $W(T,\GL_n) =
S_n$.  If $E$ is the standard representation of $\GL_n$ then the total
character of the $T-$module $E$ decomposes into a sum of linearly
independent characters $\lambda_1 + \lambda_2 + \ldots \lambda_n$ and
we get $A^*_T = \Z[t_1, \ldots t_n]$ where $t_i = c_1(\lambda_i)$.
The Weyl group $S_n$ acts on $A^*_T$ by permuting the $t_i$'s and as
result $A^*_{\GL_n} = \Z[c_1, \ldots , c_n]$ where $c_i = c_i(E)$
\cite{EGcc}.

More generally (\cite[Proposition 3.6]{eit}, \cite[Theorem 6.7]{Brion}) 
if $X$ is an algebraic space
then the restriction map $A_*^{\GL_n}(X) \to A_*^T(X)$ is an injective
homomorphism of $A^*_{\GL_n}$-modules. 
Likewise if $X$ is smooth, the restriction
map $A^*_{\GL_n} X \to A^*_T X$ is an injective homomorphism
of $A^*_{\GL_n}$-algebras. In both cases the images consist
of the elements invariant under the natural action of the Weyl group.

As a result of this discussion, we may view $A^*_{\GL_n}X$ as
a subalgebra of $A^*_T X$. In particular
we may check a formula in $A^*_{\GL_n}X$ by restricting to
$A^*_T X$. 

If $V$ is a representation of rank $r$ of
$\GL_n$
then then 
the total character of the $T-$module $V$ decomposes as sum of characters $\mu_1 + \ldots \mu_r$. Let
$l_i = c_1(\mu_i)$. 
We refer to the classes $l_1, \ldots l_r$ as the {\it Chern roots}
of $V$ and view them as elements in $A^*_T X$. Any symmetric
polynomial in the Chern roots is an element of $A^*_{\GL_n} X$.

The following easy lemma is proved for torus actions in 
\cite[Section 3.3]{eit} 
and follows in general from the projective
bundle theorem \cite[Example 8.3.4]{fu}.
\begin{lemma} \label{lem.projspace}
Let $V$ be an $r$-dimensional representation of $\GL_n$ and 
let $\Pro(V)$ be the projective space of lines in $V$. 
Then 
$$A^*_{\GL_n}(\Pro(V))= 
A^*_{\GL_n}[\xi]/(\xi^{r} + C_1 \xi^{r-1} + \ldots C_r)$$
and
$$A^*_T(\Pro(V)) = A^*_T[\xi]/\left(\prod_{i =1}^r (\xi + l_i)\right)$$
where $\xi = c_1(\Oc_{\Pro(V)}(1))$ and  $C_1, \ldots , C_r \in 
A^*_{\GL_n}$ (resp. $l_1, \ldots , l_r$)  
are the equivariant Chern classes (resp. Chern roots) 
of the representation $V$.
\end{lemma}

\section{A presentation of $\Mf^{\leq 1}_0$ as a quotient}
Let $\Mf^{\leq 1}_0$ be the stack (defined over $\Spec \Z$) of
rational curves with at most one node.
\begin{propos} \label{prop.presentation}
Let $X = \sym^2 E^* \backslash \Delta_2$
where $E$ is the defining representation of $\GL_3$, $\Delta_2$ is the
second degeneracy locus of $\sym^2 E^*$ (that is to say the locus of
quadratic forms of rank $\leq 1$).
Then the stack $\Mf^{\leq 1}_0$ is represented by the quotient stack
$$
\left[ X/ \GL_3\right]
$$ where the action of
$\GL_3$ is given by
$$
(A \cdot Q)(x)=({\det A}) Q(A^{-1} x).
$$
\end{propos}
\begin{rmk}
Proposition \ref{prop.presentation} gives another proof that
$\Mf^{\leq 1}_0$ is an algebraic stack.
\end{rmk}
Consider the natural action of $\GL_3$ on $\Pro^2 = \Proj(\Sym E^*)$.
If $L$ is a $\GL_3$-equivariant line bundle on $\Pro^2$ then
the space of global sections $H^0(\Pro^2,L)$ has a natural
$\GL_3$-module structure. If $k>0$ then
$H^0(\Pro^2, \Oc(k)) = \sym^k E^*$. The bundle $\omega^{\vee}_{\Pro^2}(-1)$
is non-equivariantly isomorphic to $\Oc_{\Pro^2}(2)$ but
$H^0(\Pro^2, \omega^{\vee}_{\Pro^2}(-1)) = \sym^2 E^* \otimes \det E$ (cf.
\cite[Exercise III.8.4]{Har}).
To prove the proposition we will show that
$$\Mf^{\leq 1}_0 = 
\left[\left(H^0(\Pro^2, \omega^{\vee}_{\Pro^2}(-1))\right) \backslash
\Delta_2\right]/\GL_3.$$

Before we prove Proposition \ref{prop.presentation} we need an easy lemma.

Let $C \xrightarrow{\pi} T$ be a rational curve with at most one node and let
$\omega_{\pi}$ be the dualizing sheaf on $C$. 
\begin{lemma}\label{lem.dual}
The ${\mathcal O}_T$-module 
$\pi_*\omega^{\vee}_{\pi}$ is locally free of rank 3
and its formation commutes with base change.
\end{lemma}
\begin{proof}[Proof of Lemma \ref{lem.dual}]
Since formation of the dualizing sheaf commutes with 
base change and $\pi \colon C \to T$ is a curve it suffices, by the
theorem of cohomology and base change, to prove 
that for every geometric point $t \to T$,
$\dim H^1(C_t, \omega_{\pi_t}^\vee) = 0$
and
$\dim H^0(C_t,\omega_{\pi_t}^\vee)= 3$.
(Here $\pi_t \colon C_t \to t$
is the restriction of $\pi \colon C \to T$.)

From Serre duality $H^1(C_t,
\omega^{\vee}_{\pi_t})=H^0(C_t, \omega^{\otimes 2}_{\pi_t})^{\vee}$.
Now, when the fiber is isomorphic to $\Pro ^1$ we have
$\omega^{\otimes 2}_{\pi_t} = \Oc(-4)$ and we do not have global
sections different from 0. When $C_t$ is singular we have that the
restriction of $\omega^{\otimes 2}_{\pi_t}$ to each of the two
components is $(\Oc(-1))^{\otimes 2}=\Oc(-2)$ and consequently the
restriction of sections to these components must be 0. Similarly we
can verify that $h^0(C_t, \omega^{\vee}_{\pi_t}) = 3$.
\end{proof}

\begin{proof}[Proof of Proposition \ref{prop.presentation}]
Let $Y$ be the category fibered in groupoids
 whose objects are pairs $(C \xrightarrow{\pi} T,
\varphi)$ where $\varphi$ is an isomorphism of $\Oc_T$ sheaves
$\varphi\colon H^0(\Pro^2, \Oc(1)) \otimes \Oc_T \cong \pi_*\omega^{\vee}_{C/T}
$ and the arrows being the obvious ones.  Each object of
$Y$ has only trivial automorphisms, so $Y$ is in fact a functor.

The action of $\GL_3$ on $H^0(\Pro^2, \Oc(1))$ induces a natural left
$\GL_3$ action on $Y$. Specifically, if $A \in
\GL_3(T)$ then
$$
A \cdot (C \xrightarrow{\pi} T, \varphi) = (C \xrightarrow{\pi} T,
\varphi \circ A^{-1})
$$

We now construct a $\GL_3$-equivariant isomorphism from $Y \to X$ (or
more precisely to the functor $\Hom(\_,X)$). Clearly $\Mf^{\leq
1}_0$ is canonically isomorphic to the quotient $[Y/\GL_3]$ so 
our proposition will follow.

Given an object $\left( C
\xrightarrow{\pi} T, \varphi \right)$ in $Y$ we have a morphism
$$
C \xrightarrow{i} \Pro(\pi_* \omega^{\vee}_{\pi})
.$$
Composing with the isomorphism $\varphi$ we get a map
$$
C \xrightarrow{i} \Pro^2_T.
$$
Next we observe that the map $i$ is an embedding. To see that it suffices
to check at geometric points of $T$. 
When the fiber $C_t \xrightarrow{\pi} t$ (at a geometric point $t
\to T$) is $\Pro^1$, we have
$$
\pi_* \omega^{\vee}_{\pi} = H^0(\Pro^1, \Oc(2))
$$
and the embedding is given by the complete linear system
$|\Oc(2)|$.

When the fiber has 2 components then the dualizing sheaf has degree
$-1$ on each component so $\omega^{\vee}_{\pi}$ is $\Oc(1)$ on each
component. The space of global sections is the 3 dimensional subspace
of sections of $\Oc(1) \oplus \Oc(1)$ that agree on the singular
point. In particular this means that the linear system $|\omega^{\vee}|$
maps each component of $C_t$ to a line in $\Pro^2$.
Moreover the images of components cannot coincide as $|\omega^{\vee}|$
separates the complement of the intersection of the two components.

Let $I \subset \Oc_{\Pro_T^2}$ be the ideal sheaf of $C$. The inclusion
$I \subset \Oc_{\Pro_T^2}$
gives a global section of $I^{\vee} = \Oc_{\Pro_T^2}(C)$ over $T$. 
By adjunction,
$$\omega_\pi^{\vee} = \omega_{\Pro_T^2}^{\vee}(-C) \otimes \Oc_C.$$
The isomorphism of global sections $$\pi_*(\omega_\pi^{\vee}) = 
H^0(\Pro^2,\Oc(1)) \otimes \Oc_T$$
gives a $\GL_3$-equivariant isomorphism 
\begin{equation} \label{eq.iso} p_*( \omega_{\Pro_T^2}^{\vee}(-C)) = 
H^0(\Pro^2, \Oc(1)) \otimes \Oc_T
\end{equation}
where $p\colon \Pro_T^2 \to T$ is the projection.

The line bundles $\omega_{\Pro_T^2}^{\vee}(-C)$ and $\Oc_{\Pro_T^2}(1)$ on
$\Pro_T^2$ both restrict to $\Oc_{\Pro^2}(1)$ on the fibers of
$p$ so a priori they differ by the pullback of a line bundle on $T$.
The equivariant isomorphism \eqref{eq.iso} of global sections then implies
that they are in fact isomorphic (with their canonical
linearizations) as $G$-line bundles on $\Pro_T^2$.
Twisting by $\Oc_{\Pro_T^2}(C) \otimes \Oc_{\Pro_T^2}(-1)$ and taking global
sections gives
a $\GL_3$-equivariant isomorphism of locally free $\Oc_T$-modules
$$p_* \omega^{\vee}_{\Pro^2_T}(-1) \to p_*(\Oc_{\Pro^2_T}(C)).$$ 
Thus the 
data $(C\to T,\varphi)$ determines a global section of
$\omega_{\Pro_T^2}^\vee(-1)$ over $T$ whose restriction to each fiber
is not in $\Delta_2$. This gives our desired $\GL_3$-equivariant map
$Y \to X$.

A $\GL_3$-equivariant inverse morphism $X \to Y$
 is given as
follows: As noted above, giving a $T$-valued point $T \to X$ is the
same as giving a global section, $s$, of
$\omega_{\Pro_T^2}^{\vee}(-1)$
 over $T$ whose restriction to each
fiber is not in $\Delta_2$. Let $C \subset \Pro_T^2$ be the
subscheme defined by the cokernel of $s^{\vee}$. 
On each fiber of $\pi$ the subscheme defined by the cokernel
of $s^{\vee}$ is subscheme of $\Pro^2_T$ cut out
by the quadratic form which is the restriction of $s$ to that fiber.
Since we assume that the restriction of $s$ is not in $\Delta_2$
it follows that $C \to
\Pro_T^2$ is a family
of at most 1-nodal rational curves of degree
$2$ such that $\Oc_{\Pro_T^2}(C) = \omega_{\Pro_T^2}^{\vee}(-1)$.
Applying adjunction again we see that $\omega_{\pi}^\vee =
\Oc_C(1)$.
 Taking global sections over $T$ we see that
$\pi_*(\omega_\pi^\vee)= \pi_* \Oc_C(1) = H^0(\Pro^2, \Oc(1)) \otimes
\Oc_T$.
\end{proof}

\section{Computation of $A^*({\Mf^{\leq 1}_0})$} \label{sec.computation}

Since $\Mf^{\leq 1}_0 = [X/\GL_3]$ where $X = \sym^2 E^* \backslash \Delta_2$
we may calculate the integral Chow ring of $\Mf^{\leq 1}_0$
as $A^*_{\GL_3} X$.

Let $\Pro^5 = \Pro(\sym^2 E^*)$. As $X$ is homogeneous for the action
of $\Gm$ with weight -1, we have that $Z:=\Pro (X)$ is a well defined
open subscheme of $\Pro^5:= \Pro \left( \sym^2 (E^*)\right)$ and there
is an induced action of $\GL_3$ on $Z$. Since the determinant acts
trivially on $\Pro^5$ the projection $p: X \to Z$ is $\GL_3$-equivariant 
and makes $X$ into the total space of the principal
$\Gm$-bundle on $Z$ associated to the line bundle
$$
{\mathcal D} \otimes \Oc(-1)
$$
where ${\mathcal D}$ is the determinant of the standard representation of 
$\GL_3$ and $\Oc(-1)$ is the tautological bundle on $\Pro^5$.

\begin{lemma} \label{lem.gmbdle}
Let $H$ denote the first Chern class of $\Oc(1)$ on $\Pro^5$. 
Then the pull-back
$$
A^*_{\GL_3}(Z) \xrightarrow{p^*} A^*_{\GL_3}(X)
$$
is surjective, and its kernel is generated by $c_1 - H$.
\end{lemma}
\begin{proof}
If $p: X \to Z$ is a $G$-equivariant
$\Gm$-bundle let ${\mathcal L}$ be the associated line bundle.
Then $X$ is the complement of the 0-section in ${\mathcal L}$ so
the basic exact sequence for equivariant Chow groups implies that
$$
A^*_G(Z)/c_1({\mathcal L}) \cong A^*_G(X).
$$
In our case, the kernel of $p^*$ is generated by $c_1\left( {\mathcal D}\otimes
\Oc(-1) \right) = c_1 -H$.

\end{proof}
\subsection{The equivariant Chow ring of $Z$}
Lemma \ref{lem.gmbdle} 
reduces the calculation of $A^*_{\GL_3}(X)$ to that of $A^*_{\GL_3}(Z)$.

Consider the following embedding
\begin{eqnarray*}
i\colon E^* &\to& \sym^2 E^*\\
\phi &\mapsto& \phi^2
\end{eqnarray*}
and its passage to the quotient
$$
i\colon \Pro(E^*) \to \Pro^5
$$
Clearly the map $i$ is $\GL_3$ equivariant and its image is $\Pro(\Delta_2)$. 
Let $j \colon Z \to \Pro^5$ be the open inclusion. Then
the basic exact sequence of $A^*{\GL_3}$-modules
$$
A^*_{\GL_3} \Pro(E^*) \xrightarrow{i_*} A^*_{\GL_3} \Pro^5 \xrightarrow{j^*} A^*_{\GL_3} Z \to 0 
$$
implies that $A^*_{\GL_3}(Z) = A^*_{\GL_3}(\Pro^5)/(\Image i_*)$.

Let $l_1, l_2, l_3$  denote the Chern roots of $E^*$, so $l_i = -t_i$.
Then the
Chern roots of $\sym^2 E^*$ are $2l_1, 2l_2, 2l_3, l_1 + l_2, l_1 +
l_3, l_2 + l_3$. By definition of $c_1, c_2, c_3$ we have
\begin{eqnarray*}
c_1 &=& -(l_1 + l_2 + l_3)\\
c_2 &=& l_1 l_2 + l_1 l_3 + l_2 l_3\\
c_3 &=& -l_1 l_2 l_3
\end{eqnarray*}

Applying Lemma \ref{lem.projspace} we see that 
$$
A^*_{\GL_3} \left( \Pro^5 \right) = \Z[c_1, c_2, c_3, H] / P(H)
$$
where
$P(H)$ is the product of linear factors
$$(H + 2l_1)(H +2l_2)(H + 2 l_3)(H + l_1 + l_2)(H + l_1 + l_3)(t +
 l_2 + l_3)$$
which can rewritten as the product
$$(H^3 - 2c_1H^2 + 4c_2H - 8c_3)(H^3 - 2c_1H^2 +
(c^2_1 + c_2)H + c_3 - c_1 c_2).$$

Let $K$ be the first chern class of $\Oc_{\Pro(E^*)}(1)$. Then
as above we have
$$
A^*_{\GL_3} \Pro(E^*) = \Z[c_1,c_2,c_3, K]/(K^3 - c_1 K^2 + c_2 K - c_3).
$$ 
In particular, $A^*_{\GL_3}(\Pro(E^*))$ is generated by the classes $1,K,K^2$ as a module over $A^*_{\GL_3}$.

\subsection{Calculation via localization}
To complete our calculation we must calculate the images
of these classes. There are a number of ways to do this. The method
we use is localization for the action of the maximal torus
$T=\Gm^3 \subset \GL_3$.

Because the restriction
map $A^*_{\GL_3} \Pro^5 \to A^*_T \Pro^5$ is injective
we can view $1, K, K^2$ as classes in $A^*_T(\Pro(E^*))$
and compute $i_*1, i_*K, i_* K^2$ in $A^*_T(\Pro^5)$. Since
$i_*$ is $\GL_3$-equivariant these images will automatically
lie in the submodule $A^*_{\GL_3}(\Pro^5)$.

Now $A^*_{T}(\Pro^5)$ is a free $A^*_T$-module generated by the
classes $1,H, \ldots H^5$. Thus any formula for $i_*K^n$ as a linear
combination of the classes $H^i$ with coefficients in $A^*_T$ can be
obtained after tensoring with $\Q$ and also localizing at the
multiplicative set of
homogeneous elements of positive degree in $A^*_T$. 
This allows us to use the explicit localization theorem 
\cite[Theorem 2]{egl}.

\begin{theorem}[\bf Explicit localization]
Let $X$ be a smooth variety equipped with an action of a torus $T$ and
let $\Qc$ be $\left( \left( A^*_T \right) ^+ \right)^{-1}A^*_T$. For every
$\alpha$ in $A^T_*(X) \otimes \Qc$ we have
$$
\alpha  = 
\sum_{F} i_{F*} \frac{i^*_F \alpha}{c^T_{\text{top}} \left( N_F X\right)}
$$
where the sum is over the components $F$ of the locus of fixed points for the action of $T$ and $i_F$ denotes the inclusion $F \to X$.
\end{theorem}

Let $\lambda_1, \lambda_2, \lambda_3$ be the basis for the
characters of $T = \Gm^3$ where $\lambda_i$ corresponds to projection to the
$i$-th factor. Then then 
the total character of the $T-$module $E^*$ decomposes as the sum of characters $$\lambda_1^{-1} + \lambda_2^{-1} + \lambda_3^{-1}.$$
Hence the Chern roots $l_1, l_2 , l_3$ of $E^*$ are $-c_1(\lambda_1), -c_1(\lambda_2), -c_2(\lambda_3)$, respectively.
Likewise then 
the total character of the $T-$module $\sym^2 E^*$ decomposes as $$\lambda_1^{-2} + \lambda_2^{-2} + \lambda_3^{-2} + \lambda_1^{-1} \lambda_2^{-1}
+ \lambda_1^{-1} \lambda_3^{-1} + \lambda_{2}^{-1} \lambda_3^{-1}.$$
Fix coordinates $[X_0 \colon X_1 \colon  X_2]$ on $\Pro^2:=\Pro(E^*)$ so
that $T$ acts by 
$$t \cdot [X_0 \colon X_1 \colon X_2]
= [\lambda_1^{-1}(t) X_0 \colon \lambda_2^{-1}(t) X_1 \colon
\lambda_3^{-1}(t)X_2]$$ and coordinates $[Z_0 \colon Z_1 \colon Z_2
\colon Z_3 \colon Z_4 \colon Z_5]$ on $\Pro^5$ so that
$T$ acts via the characters $$\lambda_1^{-2}, \lambda_2^{-2}, \lambda_3^{-2},
\lambda_{1}^{-1}\lambda_2^{-1}, \lambda_{1}^{-1}\lambda_3^{-1},
\lambda_{2}^{-1}\lambda_3^{-1}$$ on the respective coordinates.

With these actions
$\Pro^2$ has three fixed points
\begin{eqnarray*}
P_0 &:=& [1,0,0]\\
P_1 &:=& [0,1,0]\\
P_2 &:=& [0,0,1]
\end{eqnarray*}
and $\Pro^5$ has 6 fixed points $Q_0, \dots, Q_5$ defined analogously.
With our choices of coordinates, $i_*[P_j] = [Q_j]$ for $j = 0,1,2$.

Applying explicit localization to $\Pro^2$ we have
\begin{eqnarray*}
 1 &=& \sum^{2}_{k=0} i_{P_k *}\frac{i^*_{P_k} 1}{c^T_{\text{top}}(T_{P_k} \Pro^2)}\\
K &=& \sum^{2}_{k=0} i_{P_k *}\frac{i^*_{P_k} K}{c^T_{\text{top}}(T_{P_k} \Pro^2)} \\
K^2 &=& \sum^{2}_{k=0} i_{P_k *}\frac{i^*_{P_k} K^2}{c^T_{\text{top}}(T_{P_k} \Pro^2)}
\end{eqnarray*}
Let us compute the class $c^T_{\text{top}}(T_{P_0} \Pro^2)$. Local coordinates for $P_0$ are $x= X_1 / X_0, y= X_2/X_0$, these coordinates are the same for  the tangent space at $P_0$, so the action of $\Gm^3$ on  $T_{P_0} \Pro^2$ is
\begin{eqnarray*}
\Gm^3 \cdot \C^2 &\to& \C^2\\
t \cdot (x,y) &\mapsto& \left((\lambda_1\lambda_2^{-1})(t)x, 
(\lambda_1\lambda_3^{-1})(t)y\right)
\end{eqnarray*}
Consequently
$$
c^T_{\text{top}}(T_{P_0} \Pro^2) = (l_2 - l_1)(l_3 - l_1)
$$
similarly we have
\begin{eqnarray*}
c^T_{\text{top}}(T_{P_1} \Pro^2) &=& (l_1 - l_2)(l_3 - l_2)\\
c^T_{\text{top}}(T_{P_2} \Pro^2) &=& (l_1 - l_3)(l_2 - l_3)
\end{eqnarray*}

Now observe that a generator for $i^*_{P_k} \Oc(1)$ is the dual form $X_k$, so we have
\begin{eqnarray*}
i^*_{P_0}K &=& - l_1\\
i^*_{P_1}K &=& - l_2\\
i^*_{P_2}K &=& - l_3
\end{eqnarray*}

Therefore we obtain (after taking the pushforward to $\Pro^5$)
\begin{eqnarray*}
i_*1 &=& \frac{[Q_0]}{(l_2 - l_1)(l_3 - l_1)} + \frac{[Q_1]}{(l_1 - l_2)(l_3 - l_2)} + \frac{[Q_2]}{(l_1 - l_3)(l_2 - l_3)}\\
i_* K &=& \frac{- l_1[Q_0]}{(l_2 - l_1)(l_3 - l_1)} + \frac{-l_2[Q_1]}{(l_1 - l_2)(l_3 - l_2)} + \frac{-l_3 [Q_2]}{(l_1 - l_3)(l_2 - l_3)}\\
i_*K^2 &=& \frac{l^2_1[Q_0]}{(l_2 - l_1)(l_3 - l_1)} + \frac{l^2_2 [Q_1]}{(l_1 - l_2)(l_3 - l_2)} + \frac{l^2 _3 [Q_2]}{(l_1 - l_3)(l_2 - l_3)}
\end{eqnarray*}

The point $Q_0$ is the complete intersection of the  hyperplanes
cut out by the coordinate functions $Z_1, Z_2, Z_3, Z_4,Z_5$.
The $T$-equivariant fundamental classes of
these hyperplanes can be computed from the weights of the $T$-action.
For example, since $T$ acts via the character $\lambda_{2}^{-2}$ on
the coordinate $Z_1$ we see that
$[V(Z_1)] = (H + 2l_2)$ since $l_2 = c_1(\lambda_{2}^{-1})$.
Multiplying out  we obtain
\begin{eqnarray*}
 [Q_0] &=& (H +2l_2)(H +2l_3)(H + l_1 + l_2)(H + l_1 + l_3)(H + l_2 + l_3).
\end{eqnarray*}
Similar calculations show that
\begin{eqnarray*}
 \; [Q_1] &=& (H +2l_1)(H +2l_3)(H + l_1 + l_2)(H + l_1 + l_3)(H + l_2 + l_3)\\
 \; [Q_2] &=& (H +2l_1)(H +2l_2)(H + l_1 + l_2)(H + l_1 + l_3)(H + l_2 + l_3).
\end{eqnarray*}

After straight-forward computations and substituting with Chern
classes we obtain

\begin{eqnarray*}
i_*1 &=& 4(H^3 - 2c_1 H^2 + (c^2_1 + c_2)H + (c_3 - c_1 c_2))\\
i_* K &=& 2H(H^3 - 2c_1 H^2 + (c^2_1 + c_2)H + (c_3 - c_1 c_2))\\
i_*K^2 &=& H^2 (H^3 - 2c_1 H^2 + (c^2_1 + c_2)H + (c_3 - c_1 c_2))
\end{eqnarray*}

Finally we substitute $H \to c_1$ and we see that the relations in
$A^*_{\GL_3}X$ are generated by
$4c_3, 2c_1 c_3, c^2_1,c_3$.
Therefore,
$$
A^* \left( \Mf ^{\leq 1} _0 \right) = \Z [c_1, c_2, c_3] / \left( 4 c_3, 2c_1 c_3, c^2_1 c_3 \right).
$$

\section{Chow rings of the stack of reduced quadrics 
and finite extensions of $SO_n$}

Let $E$ be the defining representation of $\GL_n$ and $\Delta_{j}$ the
degeneracy locus in $\sym ^2 E^*$ of matrices with rank at most
$n-j$. Consider the following action of $\GL_n$ on  $\sym^2 E^*$:

\begin{equation} \label{eq.action}
(A \cdot Q)(x)=({\det A})^k Q(A^{-1} x)
\end{equation}
with $k \in \Z
$. 
Each subscheme $\Delta_{j}$ is invariant for this action 
and we set $X_j \colon =\sym^2 E^* \backslash \Delta_{j}$. Let
$\X_{j,k}$ denote the quotient $[X_j/\GL_n]$ where the action
of $\GL_n$ is 
given by \eqref{eq.action}.

The localization method of 
Section \ref{sec.computation} can be generalized to 
compute $A^*(\X_{n-1,k})$. Note that the stack $\X_{n-1,0}$ is the stack
of reduced quadrics in $\Pro^{n-1}$. 

On the other hand, if $j=1$ then $\X_{1,k}$ is the quotient by $\GL_n$
of the open subset $X_1$ of nondegenerate quadratic forms in $\Sym^2
E^*$. Since $X_1$ is a homogeneous space for the action of $\GL_n$
the quotient $\X_{1,k}$ is the classifying stack 
$BO(n,k)$ where $O(n,k)$ is the stabilizer of any nondegenerate quadratic form.
In our case we may identify the stabilizer of a nondegenerate quadratic form 
with 
to the closed subgroup $O(n,k) \subset \GL_n$ defined by the
condition $(\det A)^k I = A A^t$.
The groups $O(n,k)$ are extensions of $SO(n,k)$ by the cyclic 
groups $\mu_{nk-2}$ and
the techniques of \cite{Pan} can be used to calculate
the Chow ring of their classifying stacks.

\subsection{Computation of $A^*(\X_{n-1,k})$}
Let $$e_{k,n} = c_{{\rm top}}\left( (\det E)^{\otimes k} \otimes \wedge^2 E^*\right)
\in A^*_{\GL_n}$$

\begin{propos}
$A^*(\X_{n-1,k}) = \Z[c_1, \ldots , c_n]/I$
where $I$ is the ideal generated by the classes
$$\{ 2^{n-1-r}(kc_1)^{r}e_{k,n} \}_{r=0, \dots, n-1}.$$
In particular the Chow ring of the stack of reduced quadrics is
$$A^*(\X_{n,0})= \Z[c_1, \ldots , c_n]/
\left(2^{n-1}c_{{\rm top}}(\wedge ^2 E^*)\right).$$
\end{propos}

\begin{rmk} Observe if $k$ is even then $I$ is generated
by the single relation $2^{n-1}e_{k,n}$.
\end{rmk}

\begin{proof}
Let $X:= X_{n-1}$.  Since $X$ is an open set in a representation
of $\GL_n$ we can express
$A^*_{\GL_n} (X)$ as a quotient of the polynomial ring 
$\Z [c_1, \dots, c_n]$. Let $N=
\binom{n+1}{2} - 1$, $\Pro^N: = \Pro(\sym^2 E^*)$ and
$Z:=\Pro(X)$. With the action given by \eqref{eq.action} the
projection $\pi:X \to Z$ is a $\Gm$-torsor corresponding to the
line bundle $(\det E^{\otimes k}) \otimes \Oc_{\Pro^N}(-1)$.
Thus the pullback
$A^*_{\GL_n} (Z) \xrightarrow{\pi^*} A^*_{\GL_n} (X)
$ is surjective and its kernel is generated by $kc_1 - H$ where $H=
c_1(\Oc_{\Pro^N}(1))$.

Moreover we have that
$$
A^*_{\GL_n}(Z)=A^*_{\GL_n}(\Pro^N)/(\Image i_*)
$$
where $i$ is the second Veronese embedding
$$
i:\Pro(E^*) \to \Pro^N
$$
induced by
\begin{eqnarray*}
i\colon E^* &\to& \sym^2 E^*\\
\phi &\mapsto& \phi^2
\end{eqnarray*}
If we denote by $l_1, \dots, l_n$ the Chern roots of $E^*$, then the
Chern roots of $\sym^2 E^*$ are
$$
\{ 2l_i \}_{i=1, \dots, n} \cup \{ l_i + l_j \}_{1 \leq i < j \leq n}.
$$
By definition of $c_1, \dots, c_n$ we have
$$
c_k= (-1)^k s_k
$$
where $s_k$ is the $k^{\text{th}}$ symmetric polynomial in $l_1, \dots, l_n$.
Set
\begin{eqnarray*}
P(H) &=& \prod^{n}_{i=1}(H + 2l_i)\\
R(H) &=& \prod_{1 \leq i < j \leq n} (H + l_i + l_j)
\end{eqnarray*}
Because $P(H)$ and $R(H)$ are both symmetric in the $l_i$ they can be 
expressed as  polynomials with coefficients in $\Z[c_1, \dots, c_n]$. 
In particular we have
$$
P(H)= \sum_{i=0}^n (-2)^ic_i H^{n-i}. 
$$
Now applying Lemma \ref{lem.projspace} we see that
$$
A^*_{\GL_n}\left( \Pro^N \right)=\Z[c_1, \dots, c_n, H]/P(H)R(H)
$$ 
Likewise, 
if we let $K = c_1 (\Oc_{\Pro(E^*)}(1))$. 
then $A^*_{\GL_n}(\Pro(E^*))$ is generated
by the classes $1,K,K^2, \dots, K^{n-1}$.  

Let $T$ be the maximal torus in $\GL_n$ consisting of diagonal
matrices.
Applying the explicit localization formula as in Section \ref{sec.computation}
we can compute $i_* K^r \in A^*_T\Pro^N \supset A^*_{\GL_n} \Pro^N$
and we obtain the following formula
\begin{equation} \label{eq.lagrange}
i_*K^r = R(H) \sum_{j=1}^n \left( \prod_{k \neq j}\frac{H + 2l_k}{(l_k - l_j)} \right) \cdot (-l_j)^r.
\end{equation}
Viewing the sum in right-hand side of \eqref{eq.lagrange} as a polynomial of
degree $n$ in $H$ we can simplify by applying the Lagrange
Interpolation Formula. In fact the sum is the unique
polynomial of degree $n-1$ in $H$ which when evaluated at
$H=-2l_m$ (for each of $m=1, \dots, n$) equals
$(2)^{n-1}(-l_m)^r=(2)^{n-1-i}(-2l_m)^r$. Therefore,
$$
i_*K^i = 2^{n-1-r}H^{i}R(H).
$$
Moreover, the polynomial $P(H)R(H)$ is in the ideal generated by 
the $i_*K^r$'s. Consequently (after substituting $H \to kc_1$) 
we obtain that relations in $A^*_{\GL_n}X$ are generated by
$$
\{ 2^{n-1-r}(kc_1)^{r}R(kc_1)\}_{r=0, \dots, n-1}.
$$
Next, observe that
$$
R(kc_1)=\prod_{1 \leq i < j \leq n} (kc_1 + l_i + l_j)
$$
is the top Chern class of the $\GL_n$-module
$(\det E)^{\otimes k} \otimes \wedge ^2 E^*.$

In particular when $k=0$ all relations vanish, except for
the single relation $
2^{n-1}c_{top}(\wedge ^2 E^*)
$
\end{proof}

\subsection{Computation of $A^*( \X_{1,k})$}
Define polynomials $\alpha_i(H)$ in the ring $\Z[c_1, \ldots , c_n][H]$ by the formulas
\begin{equation} \label{eq.oddalpha}
 \alpha_{i}(H) = \sum^{i-1}_{j=0} \binom{n-j}{i-j}(-1)^jc_j H^{i-j} -2c_{i}
\end{equation}
when $i$ is odd
and 
\begin{equation} \label{eq.evenalpha}
 \alpha_{i}(H) = \sum^{i-1}_{j=0} \binom{n-j}{i-j}(-1)^jc_j H^{i-j}
\end{equation}
when $i$ is even.

\begin{propos}
$$A^*(\X_{1,k}) = A^*_{O(n,k)} = \Z[c_1, \ldots , c_n]/\left(\alpha_1(k c_1),
\ldots , \alpha_n(k c_1) \right).$$
\end{propos}
\begin{rmk} Observe that  $\alpha_{i}(0) = 0$ when $i$ is even and $\alpha_{i}(0) =
-2c_{i}$ when $i$ is odd, so that when $k=0$ we recover Panharipande's presentation
of $A^*_{O(n)}$.
\end{rmk}

\begin{proof}
As above, let $N=
\binom{n+1}{2} - 1$, $\Pro^N: = \Pro(\sym^2 E^*)$ and
$Z_1:=\Pro(X_1)$
With the action given by \eqref{eq.action}, the map
$X_1 \to Z_1$ is the $\Gm$-torsor corresponding to the 
line bundle $(\det E^{\otimes k}) \otimes \Oc_{\Pro^N}(-1)$.
Hence, 
$$A^*_{\GL_n}(X_1) = A^*_{\GL_n}(Z_1)/(H - k c_1).$$

We now calculate $A^*_{\GL_n}(Z_1)$ using Pandharipande's
technique. Since $Z_1 = \Pro^N \backslash \Pro(\Delta_1)$,
$A^*_{\GL_n}(Z_1) = A^*_{\GL_n}(\Pro^N)/\Image j_*$, 
where $j \colon \Pro(\Delta_1) \to  \Pro^N$ is the inclusion.

Applying \cite[Lemma 2]{Pan} with base $M$ a Chow approximation
to the classifying space $B\GL_n$ we see that 
$$
A^*_{\GL_n} (Z_1) = \Z[c_1, \dots, c_n,H]/I
$$
where $I$ is generated by $\beta_1^\prime, \dots, \beta_n^\prime$ defined as
$$
\frac{c^{\GL_n}(E^* \otimes \Oc(1))}{c^{\GL_n}(E)}= 
1 + \beta_1^\prime+ \dots + \beta_n^\prime + \dots.
$$ 
where $c^{\GL_n}$ refers to the total equivariant Chern class.
As above we can compute by restricting to torus action.
Then 
\begin{eqnarray*}
P &:=& c^T(E^* \otimes \Oc(1))= \prod^n_{i=1} (1 + H + l_i)\\
&=& \sum^n_{j=0}(-1)^j c_j(1+H)^{n-j} \\
R &:=& c^T(E)= \prod^n_{i=1}(1-l_i)\\
&=& 1 + c_1 + c_2 + \dots + c_n.
\end{eqnarray*}
where by convention $c_0=1$.

We are looking for $\beta^\prime_i$ such that
$$
P= (1 + \beta^\prime_1 + \dots + \beta^\prime_n + \dots)R,
$$

Let $P_i$ be the sum of the terms of $P$ of degree $i$.
Arguing by induction we see that 
the ideal generated by $\beta_1^\prime, \dots, \beta_n^\prime$ 
is the same as that generated by
$$
\{ \alpha_i:= P_i - c_i \}_{i=1, \dots, n}
$$
More precisely for each $i=1, \dots, n$, if $i$ is odd we have
\begin{eqnarray*}
 \alpha_{i} &=& \sum^{i-1}_{j=0} \binom{n-j}{i-j}(-1)^jc_j H^{i-j} -2c_{i}
\end{eqnarray*}
and if $i=2m$ is even we have
\begin{eqnarray*}
 \alpha_{i} &=& \sum^{i-1}_{j=0} \binom{n-j}{i-j}(-1)^jc_j H^{i-j}.
\end{eqnarray*}

Viewing the $\alpha_i$'s as polynomials in $H$, and substituting
$H \to k c_1$ we conclude
$
A^*_{\GL_n}(X_1)= \Z[c_1, \dots, c_n] / (\alpha_1(kc_1), \dots, \alpha_n(kc_1))
$
\end{proof}

\begin{rmk}
There are values of $k$ and $n$ for which 
some of the generators $ \alpha_1, \dots,
\alpha_n$ can be eliminated from the ideal of relations.  For example
if $n=3$ we have $\alpha_2 = c_1 \alpha_1$ for every $k$. On the other
hand, when $n=4$ and  $k=1$ we have
$$ (\alpha_1, \dots, \alpha_4)=(2c_1 , c^2_1, 2c_3, c_1c_3)$$
and when $k=3$
$$ (\alpha_1, \dots, \alpha_4)=(10c_1 , 5c^2_1, c^3_1 + 6c_1c_2 -2c_3, c^2_1c_2 - c_1c_3).$$
In both of these cases no generator may be eliminated.
\end{rmk}


\begin{thebibliography}{99}
\bibitem[\textbf{Br}]{Brion} M.~Brion: Equivariant Chow groups for torus
actions; Trans. Groups, \textbf{2} n. 2. 225-267 (1997).

\bibitem[\textbf{Ed-Gr1}]{EGcc} D.~Edidin, W.~Graham: Characteristic classes in the Chow ring; J. Algebraic Geom. \textbf{6} n.3, 431-443 (1997).

\bibitem[\textbf{Ed-Gr2}]{eit} D.~Edidin, W.~Graham: Equivariant intersection theory; Inv. Math. \textbf{131}, 595-634 (1998).

\bibitem[\textbf{Ed-Gr3}]{egl} D.~Edidin, W.~Graham: Localization in equivariant intersection theory and the Bott residue formula.  Amer. J. Math.  \textbf{120} no. 3, 619--636  (1998).

\bibitem[\textbf{Fulg}]{fulg} D.~Fulghesu, PhD Thesis, Scuola Normale Superiore, Pisa 2005.

\bibitem[\textbf{Ful}]{fu} W.~Fulton: \emph{Intersection theory} Springer-Verlag (Berlin), 1998.

\bibitem[\textbf{Har}]{Har} R.~Hartshorne: \emph{Algebraic Geometry} Springer-Verlag New York, 1977.

\bibitem[\textbf{Kre}]{kr} A.~Kresch: Cycle groups for Artin Stacks Inv. Math. \textbf{138} 495-536 (1999).

\bibitem[\textbf{Pan}]{Pan} R.~Pandharipande: Equivariant Chow rings of $O(k), SO(2k+1)$ and $SO(4)$; J. Reine Angew. Math. \textbf{496}, 131-148 (1998).

\bibitem[\textbf{Vis}]{vi-ap} A.~Vistoli: The Chow ring of $\Mf_2$ (Appendix to \cite{eit}); 
Inv. Math. \textbf{131} 635-644 (1998).

\end{thebibliography}
\end{document}